\documentclass[A4paper,11pt,english]{article}
\usepackage{babel}
\usepackage{graphicx}
\usepackage{amsmath}
\usepackage[utf8]{inputenc}
\usepackage{amsfonts}
\usepackage{amssymb}
\usepackage{amsmath}
\usepackage{amsthm}
\usepackage{xcolor}
\usepackage{array}
\usepackage{hyperref}
\title{\bf Convergence results in Orlicz spaces for sequences of max-product Kantorovich sampling operators}
\author{ {\bf Lorenzo Boccali} \hskip1cm {\bf Danilo Costarelli} \hskip1cm {\bf Gianluca Vinti} \\ 
	Department of Mathematics and Computer Science \\
	University of Perugia\\
	1, Via Vanvitelli, 06123 Perugia, Italy \\ 
	{\small {\tt lorenzo.boccali@unipg.it}} - {\small {\tt danilo.costarelli@unipg.it}} - {\small {\tt gianluca.vinti@unipg.it}} }
\date{}
\newtheorem{prop}{Proposition}[section]
\newtheorem{definizione}{Definition}[section]
\newtheorem{lemma}{Lemma}[section]
\newtheorem{teorema}{Theorem}[section]
\newtheorem{cor}{Corollary}[section]
\theoremstyle{definition}
\newtheorem{remark}{Remark}[section]
\newcommand{\N}{\mathbb{N}}
\newcommand{\R}{\mathbb{R}}
\newcommand{\Z}{\mathbb{Z}}
\newcommand{\la}{\lambda}
\newcommand{\fhi}{\varphi}
\newcommand{\ep}{\varepsilon}
\newcommand{\Om}{\Omega}
\newcommand{\assolutol}{\lvert}
\newcommand{\assolutor}{\rvert}
\newcommand{\norma}{\|}

\newcommand{\Lp}{L_{+}^{\varphi}(\Omega)}
\newcommand{\I}{I^{\varphi}}
\newcommand{\op}{K_{n}^{\chi}}
\newcommand{\mom}{m_{0}(\chi)}

\newcommand{\integrale}{\int_{\Omega}}
\newcommand{\integralek}{\int_{k/n}^{(k+1)/n}}
\newcommand{\integraleab}{\int_{a}^{b}}
\newcommand{\integraler}{\int_{\R}}
\newcommand{\integralem}{\int_{\assolutol x \assolutor>M}}
\newcommand{\integraleb}{\int_{B}}
\newcommand{\V}{\bigvee}
\newcommand{\Vz}{\bigvee_{k \in \mathbb{Z}}}
\newcommand{\Vj}{\bigvee_{k \in \mathcal{J}_{n}}}
\newcommand{\Vg}{\bigvee_{k \in \mathcal{J}_{n}^{\gamma }}}
\newcommand{\J}{\mathcal{J}}
\newcommand{\mb}{m_{\beta}(\chi)}
\newcommand{\rapporto}{\dfrac{\lambda}{a_{\chi}}}
\begin{document}
\maketitle
\begin{abstract}
\noindent In this paper, we provide a unifying theory concerning the convergence properties of the so-called max-product Kantorovich sampling operators based upon generalized kernels in the setting of Orlicz spaces. The approximation of functions defined on both bounded intervals and on the whole real axis has been considered. Here, under suitable assumptions on the kernels, considered in order to define the operators, we are able to establish a  modular convergence theorem for these sampling-type operators. As a direct consequence of the main theorem of this paper, we obtain that the involved operators can be successfully used for approximation processes in a wide variety of functional spaces, including the well-known interpolation and exponential spaces. This makes the Kantorovich variant of max-product sampling operators suitable for reconstructing not necessarily continuous functions (signals) belonging to a wide range of functional spaces. Finally, several examples of Orlicz spaces and of kernels for which the above theory can be applied are presented.     
\end{abstract}
\medskip\noindent
{\small {\bf AMS subject classification:} 41A25, 41A35 \newline
{\small {\bf Key Words:} Max-product Kantorovich sampling operators; generalized kernels; Orlicz spaces; modular convergence; Luxemburg norm.}
\section{Introduction}
In view of their many applications in signal and image processing, sampling operators represent one of the most studied topics in Approximation Theory. Beginning in the first half of the $1900$s with the classical sampling theorem (see, e.g., \cite{shannon,Butzer}), with the major contributions of Whittaker, Kotelnikov and Shannon, their possible use for signal reconstruction or approximation have attracted strong interest in the scientific community. Since the required assumptions on the signal to be approximated limit its applicability to solving real-world problems, possible approximate versions of the above theorem have been studied since the $1960$s. A first rigorous theory in this direction has been developed by the German mathematician P. L. Butzer and his school in Aachen with the introduction of the so-called generalized sampling operators (series) which are characterized by good approximating properties (see, e.g., \cite{ries,Butzer2,Butzer3,Butzer4,Butzer5,kivinukk}). In the following years, many authors have been concerned with the study of alternative versions of the above operators acting in a wide variety of functional spaces. Recently, in many papers (see, e.g., \cite{2.5}-\cite{4.4}), the max-product approach applied to several approximation operators has been studied by Coroianu and Gal. More in details, by replacing the series (or sum in case of finite terms), defining a family of discrete approximation operators, with the symbol $\bigvee$, which denotes the supremum (or maximum, in the finite case) of a given set of real numbers, it is possible to obtain the nonlinear (more precisely, sub-additive) version of any family of linear operators.  In addition to preserving the approximation capability of the corresponding linear counterparts, max-product operators often return better orders of approximation (for more details, see, e.g., \cite{Karakus,01,CDGV} or Remark \ref{remark5.1} of Section \ref{esempigrafici}). For this reason the max-product method is widely used for applications to various theoretical fields of Approximation Theory, such as function approximation (see, e.g., \cite{006,Duman,7,7.8,6}), theory of neural networks (see, e.g., \cite{5.1,5.2}) and also fuzzy sets theory (see, e.g., \cite{04,9}). Within the framework of sampling theory, in \cite{2}, the authors introduced  the max-product version of Kantorovich sampling operators based on generalized kernels for the approximation of univariate signals. We recall here that the linear version of the Kantorovich sampling operators (series) have been introduced and studied by Bardaro et al. \cite{004} in univariate form, where approximation results in Orlicz spaces have been established. In the case of a uniform sampling-scheme, these series are defined (in the univariate setting) as follows:
\begin{equation*}
S_{w}^{\chi}\bigl(f\bigr)(x) := \sum_{k \in \Z}\chi(wx-k) \left[ w \int_{k/w}^{(k+1)/w} f(t) \ dt\right], \quad x \in \R, \quad w >0,
\end{equation*}
where $f: \R \rightarrow \R$ is a locally integrable function (signal) such that the above series is convergent for every $x \in \R$. Here, $\chi: \R \rightarrow \R$ is a kernel function satisfying the usual conditions of the discrete approximate identities (see, e.g., \cite{02,7.9}). It is well known that, in view of their averaged nature, sampling operators of Kantorovich-type represent a useful tool to obtain an approximate version of the celebrated Whittaker-Kotelnikov-Shannon sampling theorem, proving particularly suitable for the reconstruction of signals that are not necessarily continuous (see, e.g., \cite{004,Orlova,4.6}). In general, Kantorovich-variants of sampling operators have been extensively studied in the last fifteen years from both a theoretical and applied point of view, returning a wide variety of applications to image processing in several fields, such as biomedical and engineering (see, e.g., \cite{cagini}).  \newline 
As mentioned above, Kantorovich sampling operators of max-product type have been introduced by Coroianu et al. \cite{2} and they are formally expressed by:
\begin{equation*}
\op \bigl(f\bigr)(x):=\dfrac{\displaystyle \bigvee_{k \in \Z}\chi(nx-k) \left[n \integralek f(t) \ dt\right]}{\displaystyle \bigvee_{k \in \Z} \chi(nx-k)}, \quad x \in \R, \quad (\textnormal{I})
\end{equation*}    
where $f: \R \rightarrow \R$ is a locally integrable function and $\chi: \R \rightarrow \R$ is a generalized kernel satisfying suitable assumptions. By adopting the so-called moment-type approach in its max-product variant, obtained by replacing the sum (or series) in the usual definition of the discrete absolute moments of a given kernel (see, e.g., \cite{ries}), it has been possible to investigate the approximation properties of max-product Kantorovich sampling operators $(\textnormal{I})$ in several functional spaces, including the $L^{p}$-spaces, $1 \le p < +\infty$, and suitable spaces of continuous functions. In both cases, qualitative and quantitative convergence properties have been obtained. In order to prove these results, unlike the corresponding linear counterpart, it is not necessary to assume that $\chi$ is an approximate identity. Indeed, the theory developed  in \cite{2} still holds if we consider measurable and bounded kernels with in addition a finiteness assumption on the so-called generalized absolute moment for a given order $\beta >0$. 
\newline The goal of this paper is to show that, keeping these assumptions on $\chi$, we can extend the results proved in \cite{4.5,coroianu1,2} in a twofold sense, on the one hand considering kernels not necessarily having compact support or generated by sigmoidal functions (as, e.g., in \cite{2.5}) and on the other hand approximating functions which belong to the Orlicz space $L^{\fhi}$ generated by a convex $\fhi$-function. As it is well-known from the literature (see, e.g., \cite{7.3,7.1}), these general spaces include the usual $L^{p}$-spaces and several other well-known functional spaces, such as the interpolation and the exponential spaces. In particular, we will prove a converge result for sequences of max-product Kantorovich sampling operators $(\textnormal{I})$ with respect to the most natural notion of convergence which can be defined in $L^{\fhi}$, i.e., the so-called modular convergence. To reach this, after recalling in Section \ref{preliminari} the definition of the involved operators and all the auxiliary results used in this paper, we will establish a modular inequality for the same family of non-linear operators. Next, as often happens in the case of approximation processes in Orlicz spaces (whose theory is briefly recalled in Section \ref{orlicz}), we will need to investigate the modular convergence for continuous functions on compact intervals (or continuous functions with compact support if $f$ is defined on $\R$) so that a density approach can be applied. In this way, the main theorem of the present paper given in Section \ref{risultaticonvergenza} easily follows. In Section \ref{esempigrafici}, we will consider some particular cases of Orlicz spaces, finding again the convergence results with respect to the $L^{p}$-norm which have been proved in \cite{2}. Then, we will provide several examples of kernels for which the above theory can be applied, together with the corresponding graphical representations.  
\section{Preliminary results and notations}
\label{preliminari}
From now on, in the whole paper, by the notation $f: \Om \rightarrow \R$ we will denote real functions defined on $\Om \subseteq \R$. If $\Om$ is a compact interval $[a,b]$, then $C(\Om)$ denotes the set of all continuous functions defined on $\Om$, while in the case of $\Om=\R$, $C(\Om)$ denotes the set of all uniformly continuous and bounded functions defined on $\Om$. In both cases, $C(\Om)$ is a normed linear space with the usual sup-norm, that we will denote by $\norma \cdot \norma_{\infty}$. Moreover, let us denote by $C_{+}(\Om)$ the subspace of $C(\Om)$ of the non-negative functions and by $C_{c}^{+}(\Om) \subset C_{+}(\Om)$ the subspace whose elements have compact support.      
\newline Now, we recall that the symbol $\V$, in the literature (see, e.g., $\cite{3,4,01,6,7,Shen}$), denotes the supremum of any set of real numbers $\{A_{k}: k \in \J\}$ with respect to a given $\J \subseteq \Z$ set of indexes, i.e., \[\V_{k \in \mathcal{J}}A_{k}:= \sup \ \{A_{k} \in \R, k \in \J\}.\]
Obviously, if $\mathcal{J}$ has finite cardinality, the supremum reduces to a maximum.
As we will see later, the above definition is crucial in order to define the kind of sampling operators that we will study in this paper in the general context of Orlicz spaces. Such operators have been firstly introduced and studied in \cite{2}, where the approximation of (non-negative) continuous and $L^{p}$-integrable functions has been considered.
\newline In order to recall the definition of the above-mentioned operators, we firstly give the notion of generalized kernels, i.e., any bounded and measurable function $\chi: \R \rightarrow \R$ such that the following properties are satisfied: 
\newline \newline  $(\chi_{1})$ the generalized absolute moment of order $\beta$ of $\chi$ is finite for some $\beta > 0$, i.e.,
\[\mb:= \sup_{x \in \R}\Vz \left| \chi(x-k) \right| \cdot \left |x-k\right| ^{\beta}< +\infty;\] 
\newline $(\chi_{2})$ if $\Om=[a,b]$, then we have:
\[\inf_{x \in [-3/2,3/2]} \chi(x) =: a_{\chi}>0;\]
or \newline $(\chi_{2}')$ in the case of $\Om= \R$, we have: \[\inf_{x \in [-1/2,1/2]} \chi(x) =: a_{\chi}>0.\]
\newline Conditions $(\chi_{1})$ and $(\chi_{2}')$ have been firstly introduced in \cite{1} in order to establish a pointwise and uniform convergence result and quantitative estimates for generalized sampling operators of max-product type. Similarly to above, as it was shown in \cite{2}, keeping the same assumptions on $\chi$ generates similar convergence properties for the max-product sampling operators of Kantorovich-type, with respect to both the norm $\norma \cdot \norma_{\infty}$ and the usual $L^{p}$-norm $\norma \cdot \norma_{p}$, $1 \le p < +\infty$, except for the bounded case $\Om=[a,b]$, where it is necessary to replace  condition $(\chi_{2}')$ with the slightly stronger assumption $(\chi_{2})$.
\newline In Section \ref{risultaticonvergenza}, we will prove that, considering kernels of the same type, it is possible to provide an extension of the results obtained in \cite{2} in the more general framework of Orlicz spaces generated by convex $\fhi$-functions, which include the $L^{p}$-spaces as some special instances.     
\newline We recall here some useful results concerning the kernel $\chi$ in the max-product setting.
\begin{lemma} [\cite{1}, \textnormal{Lemma 2.1}]
\label{lemma1}
Let $\chi: \R \rightarrow \R$ be a bounded function such that $\chi(x) = \mathcal{O}(\assolutol x \assolutor^{-\alpha})$, as $\assolutol x \assolutor \rightarrow +\infty$, for some $\alpha >0$. Then it turns out that:
\[m_{\beta}(\chi) < +\infty,\]
for every $0 \le \beta \le \alpha$. 
\end{lemma}
\begin{lemma} [\cite{1}, \textnormal{Lemma 2.2}]
\label{lemmamomenti}
Let $\chi: \R \rightarrow \R$ be a bounded function which satisfies $(\chi_{1})$ for some $\beta>0$. Then there holds: \[m_{v}(\chi) < +\infty,\]
for every $0 \le v \le \beta$. In particular, we have $ m_{0}(\chi)\le \norma \chi \norma_{\infty}$. 
\end{lemma} 
\begin{lemma}
\label{Lemma1}
Let $ \chi: \R \rightarrow \R$ be a given function which satisfies condition $(\chi_{2}')$. Then the following inequality:
\begin{equation}
\label{sup}
\V_{k \in \Z} \chi (nx-k) \ge a_{\chi} > 0, 
\end{equation}
holds for all $x \in \R$ and for every $n \in \N$. \newline In addition, let $[a,b]$ be a compact interval. If $\chi$ satisfies assumption $(\chi_{2})$, then for all $x \in [a,b]$ there holds:
\begin{equation}
\label{max}
\V_{k \in \J_{n}} \chi(nx-k) \ge a_{\chi} >0,
\end{equation}
for every $n \in \N$ sufficiently large, where $\J_{n} := \{k \in \Z: \lceil{na}\rceil \le k \le \lfloor{nb}\rfloor -1\}$. Here, $\lceil{\cdot}\rceil$ and $\lfloor{\cdot}\rfloor$ denote respectively, the ``ceiling" and the ``integral part" of a given number. Finally, $a_{\chi}$ is the constant from condition $(\chi_{2}')$ and $(\chi_{2})$, respectively. 
\end{lemma}
The proof of this lemma is similar to that one of Lemma 2.3 in \cite{1}.
\begin{remark}
As noted in \cite{2} too, it is possible to weaken the assumptions on the kernel $\chi$ by directly assuming that both the inequalities (\ref{sup}) and (\ref{max}) hold, together with condition $(\chi_{1})$, instead of conditions $(\chi_{2})$ and $(\chi_{2}')$. Indeed, in this setting, the approximation results which we will prove in this paper still hold. It follows that we can also consider kernels which are not necessarily  bounded from below by a strictly positive constant on $[-1/2,1/2]$ or $[-3/2,3/2]$.   
\end{remark}
We are now able to recall the definition of the operators considered in this paper. In the following definition $\J_{n} = \Z$, if $\Om=\R$ and $\J_{n}=\{k \in \Z : \lceil{na}\rceil \le k \le \lfloor{nb}\rfloor -1\}$, if $\Om=[a,b]$. 
\begin{definizione}
Let $\chi: \R \rightarrow \R$ be a fixed generalized kernel. We define the family of max-product Kantorovich sampling operators $(\op)_{n \in \N}$ based upon $\chi$ as follows:
\[\op (f) (x):= \dfrac{\displaystyle \Vj \chi(nx-k) \left[n \displaystyle \int_{k/n}^{(k+1)/n}f(t) \ dt\right]}{\displaystyle \Vj \chi(nx-k)}, \quad x \in \Om,\]
where $f : \Om \rightarrow \R$ is a locally integrable function.  
\end{definizione}
From now on, in the whole paper, we will always assume that the kernel $\chi$ satisfies $(\chi_{1})$ and $(\chi_{2})$ or $(\chi_{2}')$ for $\Om=[a,b]$ and $\Om=\R$, respectively. Therefore, since in this setting  Lemma \ref{Lemma1} holds, we always have $\Vj \chi(nx-k) > 0$, for all $x \in \Om$. Moreover, it is easy to see that, $\op (f)$ is well-defined, for example, if $f$ is bounded. Indeed, in this case, by Lemma \ref{lemmamomenti}, we have:
\begin{equation}
\left| \op (f) (x) \right| \le \dfrac{\norma f \norma_{\infty}}{a_{\chi}} \mom < +\infty,
\end{equation}
for every $x \in \Om$.
\newline We summarize some important proprieties for $\op$ in the following lemma.
\begin{lemma} [\cite{2}, \textnormal{Lemma 2.6}]
\label{lemmaproprietà}
Let $\chi$ be a fixed generalized kernel. Let $f$, $g: \Om \rightarrow [0, +\infty)$ be non-negative and bounded functions on $\Om$. The following properties hold for all $n \in \N$:
\newline 
\newline (i) if $ f \le g$, then $\op (f) \le \op (g)$;
\newline (ii) $\op(f+g) \le \op (f) + \op (g)$, i.e., $\op$ is sub-additive;
\newline (iii) $|\op (f) - \op (g)| \le \op (|f-g|)$;
\newline (iv) $\op (\lambda f)= \lambda \op (f)$ for each $\lambda \ge 0$, i.e., $\op$ is positive homogeneous. 
\end{lemma}
Let us now recall the following convergence theorem for the operators $(\op)_{n \in \N}$ in the classical setting, i.e., pointwise and uniform convergence in suitable spaces of continuous functions. 
\begin{teorema} [\cite{2}, \textnormal{Theorem 3.1}]
\label{convergenzauniforme}
Let $\chi$ be a given generalized kernel. Let $f: \Om \rightarrow [0, +\infty)$ be a non-negative and bounded function, continuous at $x \in \Om$. Then
\[\lim_{n \rightarrow +\infty} \op(f)(x)=f(x).\]
In addition, if $f \in C_{+}(\Om)$, then 
\[\lim_{n \rightarrow +\infty}\|\op(f)-f\|_{\infty}=0.\] 
\end{teorema}
\section{The Orlicz Spaces}
\label{orlicz}
In this section, we recall some basic notions concerning Orlicz spaces, which have been introduced by the Polish mathematician W. Orlicz as a natural extension of Lebesgue spaces (see, e.g., \cite{7.1,7.2,002}). Let us denote by $M(\Om)$ the set of all (Lebesgue-)measurable real functions defined on $\Om$, where $\Om$ can be (again) a bounded interval or the whole real axis.  \newline In what follows, we consider a function $\fhi:\R_{0}^{+} \rightarrow \R_{0}^{+}$ which satisfies the following assumptions:
\newline
\newline $(\fhi1)$ $\fhi(0)=0$, $\fhi(u) >0$ for every $u >0$;
\newline $(\fhi2)$ $\fhi$ is continuous and non-decreasing on $\R_{0}^{+}$;
\newline $(\fhi3)$ $\displaystyle \lim_{u \rightarrow +\infty} \fhi(u) = +\infty.$
\newline One can call this function as $\fhi$-function. For a fixed $\fhi$-function $\fhi$, it is possible to define the functional $\I: M(\Om) \rightarrow [0, +\infty]$ as follows:
\[\I[f]:=\integrale \fhi\bigl(\assolutol f(x) \assolutor\bigr) \ dx, \quad f \in M(\Om).\]   
As it is well known, $\I$ is a modular functional on $M(\Om)$. Now, we can introduce the Orlicz space generated by $\fhi$ as the set:
\[L^{\fhi}(\Om)=\{f \in M(\Om): \exists \ \la >0 \ such \ that \ \I[\la f] < +\infty\}.\]
From the literature, if $\fhi$ is convex, it turns out that $L^{\fhi}(\Om)$ is a normed linear space with norm $\norma \cdot \norma_{\fhi}$, that is known as Luxemburg norm, and it is defined by:
\[\norma f \norma_{\fhi}:= \inf\left\{\la >0: \I[f/\la] \le 1\right\}, \quad f \in L^{\fhi}(\Om). \]
An important and useful subspace of $L^{\fhi}(\Om)$ is
\[E^{\fhi}(\Om):= \{f \in M(\Om): \I[\la f] < +\infty, \ for \ every \ \la >0\},\]
i.e., the so-called space of all finite elements of $L^{\fhi}(\Om)$. In general, $E^{\fhi}(\Om)$ is a proper subspace of $L^{\fhi}(\Om)$ but these two spaces coincide if and only if the inequality:
\[\fhi(2u) \le M \fhi(u), \quad u \in \R_{0}^{+},\]
holds for some $M>0$. This condition on $\fhi$ is known in the literature as the $\Delta_{2}$-condition. 
\newline Examples of $\fhi$-functions for which the $\Delta_{2}$-condition is satisfied are $\fhi(u)=u^{p}$, $1 \le p < +\infty$, or $\fhi_{\alpha,\beta}(u)=u^{\alpha}\log^{\beta}(u+e)$, for $\alpha \ge 1$ and $\beta >0$, which generate the usual $L^{p}$-spaces and the well-known interpolation spaces $L^{\alpha} \log^{\beta}L(\Om)$, respectively. On the other hand, it is easy to prove that, for example, for $\fhi_{\gamma}(u)=e^{u^{\gamma}}-1$, $\gamma >0$, the $\Delta_{2}$-condition is not fulfilled. Hence, each of the corresponding Orlicz spaces $L^{\fhi_{\gamma}}(\Om)$, known as the exponential spaces, properly contains its subspace of the finite elements $E^{\fhi_{\gamma}}(\Om)$. 
\newline Now, we recall that it is possible to introduce in $L^{\fhi}(\Om)$ a natural notion of convergence based on the definition of the modular functional $\I$, i.e., the so-called modular convergence. From now on, we will say that a net of functions $(f_{w})_{w>0} \subset L^{\fhi} (\Om)$ is modularly convergent to a function $f \in L^{\fhi}(\Om)$ if:
\[\lim_{w \rightarrow +\infty} \I\bigl[\la\bigl(f_{w}-f\bigr)\bigr]=0,\]
for some $\la >0$. Moreover, the Luxemburg norm $\norma \cdot \norma_{\fhi}$ defines itself a notion of convergence in $L^{\fhi}(\Om)$ which is called the (Luxemburg) norm-convergence. Namely, a net of functions $(f_{w})_{w>0} \subset L^{\fhi}(\Om)$ is said norm convergent to a function $f \in L^{\fhi}(\Om)$, i.e., $\norma f_{w} - f\norma_{\fhi} \rightarrow 0$ for $ w \rightarrow +\infty$, if:
\[\lim_{w \rightarrow +\infty} \I\bigl[\la\bigl(f_{w}-f\bigr)\bigr]=0,\]
for every $\la >0$. From the above definitions, it is clear that the Luxemburg norm convergence is stronger than the modular convergence. Further, it is possible to prove that, assuming that $\fhi$ satisfies the $\Delta_{2}$-condition, these two notions of convergence are equivalent. Moreover, the vice versa holds.  
\newline Finally, in order to reach the purpose of this paper, we need to recall the following density properties of the Orlicz spaces. Namely, it is easy to prove that $C([a,b])\subset L^{\fhi}([a,b])$ and $C_{c}(\R) \subset L^{\fhi}(\R)$, and moreover, it turns out that $C([a,b])$ and $C_{c}(\R)$  are dense in $L^{\fhi}([a,b])$ and $L^{\fhi}(\R)$, respectively, with respect to the topology induced by the modular convergence. Similarly to above, we can prove that $C_{+}([a,b])$ and $C_{c}^{+}(\R)$ are modularly dense in $L^{\fhi}_{+}(\Om)$, with $\Om=[a,b]$ or $\Om=\R$, respectively, where $L^{\fhi}_{+}(\Om)$ denotes the subspace of $L^{\fhi}(\Om)$ of the non-negative functions. 
\newline For further results and details concerning the Orlicz spaces, one can see, e.g., \cite{7.3,7.5,7.4,7.6,7.7}. 
Finally, before to consider the operators $(\op)_{n \in \N}$ in the Orlicz spaces, we need to establish the following proposition in which the max-product symbol $\V$ is related to the convex $\fhi$-function $\fhi$.
\begin{prop}
\label{proposizione}
Let $\fhi$ be a convex $\fhi$-function. Then the following inequality:
\begin{equation*}
\fhi\left(\bigvee_{k \in \J} A_{k}\right) \le \bigvee_{k \in \J} \fhi\bigl(2 A_{k}\bigr),
\end{equation*} 
holds for any set of indexes $\J \subseteq \Z$, and $A_{k} \ge 0$, $k \in \J$.
\end{prop}
\begin{proof}
If $\mathcal{J}$ is finite, it is easy to observe that the following equality:
\begin{equation*}
\fhi\left(\bigvee_{k \in \J} A_{k}\right) = \bigvee_{k \in \J} \fhi\bigl(A_{k}\bigr), 
\end{equation*}
holds, since $\fhi$ is non-decreasing. Thus, the above inequality is obvious. Now, let us consider the case $\mathcal{J}=\Z$ (if $\mathcal{J}$ is any infinite set, the proof is similar). Let $\varepsilon >0$ be fixed. From the definition of the symbol $\V$, it follows that there exists $\bar{k} \in \Z$ such that $A_{\bar{k}}> \bigvee_{k \in \Z} A_{k} - \varepsilon$. Therefore, by the properties of $\fhi$, we can write what follows:
\begin{equation*}
\begin{split}
\fhi\left(\bigvee_{k \in \Z} A_{k} \right)  & \le \fhi\bigl(A_{\bar{k}}+\varepsilon\bigr)=\fhi\left(\frac{1}{2}2A_{\bar{k}}+\frac{1}{2}2\varepsilon\right) \le \frac{1}{2}\bigl\{\fhi\bigl(2A_{\bar{k}}\bigr)+\fhi\bigl(2\varepsilon\bigr)\bigr\} \\
 & \le \fhi\bigl(2A_{\bar{k}})+\varepsilon \fhi(2) \le \bigvee_{k \in \Z} \fhi\bigl(2A_{k}\bigr)+ \varepsilon \fhi(2).
\end{split}
\end{equation*}
Thus, the assertion follows by the arbitrariness of $\varepsilon>0$.
\end{proof}   
\section{Convergence results in Orlicz spaces}
\label{risultaticonvergenza}
In this section, we study the convergence proprieties of the max-product Kantorovich sampling operators in the general setting of Orlicz spaces. 
\newline First of all, we need to prove an inequality for these non-linear operators with respect to the modular functional $I^{\varphi}$. We will see that such inequality is crucial in order to prove the desired modular convergence result for the family $(K_{n}^{\chi})_{n \in \N}$ in $L^{\fhi}_{+}(\Om)$, with $\Om =[a,b]$ or $\Om=\R$.
\newline In all that follows in this section, if otherwise not stated, the $\varphi$-function $\varphi$ will be convex and the kernel $\chi$ from the definition of $K_{n}^{\chi}$ will belong to $L^{1}(\R)$.
\begin{teorema}
\label{disuguaglianzamod}
For every $f$, $g \in \Lp$ and $\lambda > 0$, it turns out that:
\[\I\bigl[\lambda\bigl(\op (f)- \op (g)\bigr)\bigr] \le \dfrac{\|\chi\|_{1}}{m_{0}(\chi)}\I\left[\dfrac{\mom}{a_{\chi}}2\lambda\bigl(f-g\bigr)\right],\]
for $n \in \N$ sufficiently large. 
\end{teorema}
\begin{proof}
By the property $(iii)$ of the operator $\op$ established in Lemma \ref{lemmaproprietà} and the inequalities of Lemma \ref{Lemma1}, we can write what follows for sufficiently large $n \in \N$:
\begin{equation*}
\begin{split}
\I\bigl[\lambda\bigl(\op(f)-\op(g)\bigr)\bigr]&=\integrale \varphi\bigl(\lambda\assolutol \op(f)(x)-\op(g)(x)\assolutor\bigr)  \ dx \le \integrale \varphi\bigl(\lambda\op(\assolutol f-g \assolutor)(x)\bigr)  \ dx \\
& \le \integrale \varphi\left(\dfrac{\lambda}{a_{\chi}}\Vj \assolutol\chi(nx-k) \assolutor \left[n \integralek \assolutol f(t)-g(t) \assolutor \ dt \right] \right) \ dx
\end{split}
\end{equation*} 
Now, from Proposition \ref{proposizione}, we immediately get that:
\begin{equation}
\label{disuguaglianzaphi}
\varphi\left(\V_{k \in \J}A_{k}\right) \le \V_{k \in \J} \varphi\bigl(2A_{k}\bigr),
\end{equation}
for any set of indexes $\J$. Thus, by using (\ref{disuguaglianzaphi}), the convexity of $\varphi$ with the fact that $\left| \chi(nx-k)\right| \le \mom$, for every $k \in \J_{n}$, and the Jensen inequality (see, e.g., \cite{5}), we obtain: 
\begin{equation*}
\begin{split}
\I\bigl[\lambda\bigl(\op(f)-\op(g)\bigr)\bigr] &\le \integrale \Vj \varphi\left(2\dfrac{\lambda}{a_{\chi}} \assolutol \chi(nx-k) \assolutor \left[ n \integralek \assolutol f(t)-g(t) \assolutor \ dt \right] \right) \ dx \\
& \le \integrale \Vj \dfrac{\assolutol \chi(nx-k) \assolutor}{\mom}  \ \varphi\left(2\dfrac{\mom}{a_{\chi}} \lambda \left[n \integralek \assolutol f(t)-g(t) \assolutor \ dt\right]\right) \ dx \\
&\le  \integrale dx \left\{\Vj \dfrac{\assolutol \chi(nx-k) \assolutor}{\mom} \left[ n \integralek \varphi\left(2\dfrac{\mom}{a_{\chi}} \lambda \assolutol f(t)-g(t) \assolutor \right) \ dt\right] \right\} \\ 
& \le (\mom)^{-1} \integrale \left[ n \sum_{k \in \J_{n}} \assolutol \chi(nx-k) \assolutor \integralek \fhi \left(2\dfrac{\mom}{a_{\chi}} \la \assolutol f(t)-g(t) \assolutor \right) \ dt \right] \ dx \\
\end{split}
\end{equation*}
Now, it is easy to observe that in the case $\Om=\R$ we can use the Fubini-Tonelli theorem in order to interchange the integral with the series in the above computations. Thus, we get:
\begin{equation}
\label{disuguaglianzamodulare}
\begin{split}
\I\bigl[\la\bigl(\op(f)-\op(g)\bigr)\bigr] & \le (\mom)^{-1} \integrale \left[ n \sum_{k \in \J_{n}} \assolutol \chi(nx-k) \assolutor \integralek \fhi \left(2\dfrac{\mom}{a_{\chi}} \la \assolutol f(t)-g(t) \assolutor \right) \ dt \right] \ dx \\
& = (\mom)^{-1} \sum_{k \in \J_{n}} \integralek \fhi\left(2\dfrac{\mom}{a_{\chi}}\la \assolutol f(t) -g(t) \assolutor\right) \ dt \integrale n \assolutol \chi(nx-k) \assolutor \ dx
\end{split}
\end{equation}
Moreover, for every $k \in \J_{n}$, by using the change of variable $y=nx-k$, we have:
\[ \integrale n  \assolutol \chi(nx-k) \assolutor \ dx \le \int_{\R} \assolutol \chi(y) \assolutor \ dy = \|\chi\|_{1}.\]
Hence, using this result in (\ref{disuguaglianzamodulare}), we finally obtain:
\begin{equation*}
\begin{split}
\I\bigl[\lambda\bigl(\op(f)-\op(g)\bigr)\bigr] & \le (\mom)^{-1} \norma \chi \norma_{1} \sum_{k \in \J_{n}} \integralek \fhi\left(2\dfrac{\mom}{a_{\chi}} \la \assolutol f(t)-g(t) \assolutor\right) \ dt \\
& \le  \dfrac{\norma \chi \norma_{1}}{\mom} \integrale \fhi\left(2\dfrac{\mom}{a_{\chi}}\la \assolutol f(t)-g(t) \assolutor\right) \ dt \\
& = \dfrac{\|\chi\|_{1}}{\mom} \I\left[\dfrac{\mom}{a_{\chi}} 2\la (f-g)\right],
\end{split}
\end{equation*}
for every $n \in \N$ sufficiently large. This completes the proof.
\end{proof}   
Now, in order to establish a modular convergence theorem for the max-product Kantorovich sampling operators in $L^{\fhi}_{+}(\Om)$, we need to test the (Luxemburg) norm-convergence in case of functions belonging to $C_{+}(\Om)$, when $\Om=[a,b]$, and in case of functions $f \in C_{c}^{+}(\Om)$, when $\Om=\R$. Hence, we can give the proof of the following.
\begin{teorema}
\label{convergenzaforte}
Let $f \in C_{+}(\Om)$, $\Om = [a,b]$, be fixed. Then for every $\lambda > 0$, it turns out that:
\[\lim_{n \rightarrow +\infty} \I \bigl[\lambda\bigl(\op(f)-f\bigr)\bigr]=0.\] 
In addition, in case of $f \in C_{c}^{+}(\Om)$, with $\Om=\R$, we also have:
\[\lim_{n \rightarrow +\infty} \I \bigl[\lambda\bigl(\op(f)-f\bigr)\bigr]=0,\]
for every $\lambda >0$.
\end{teorema}
\begin{proof}
Let us begin by considering the case $\Om=[a,b]$. First of all, it is easy to see that $\op(f)$ belongs to $E^{\fhi}([a,b]) \subseteq L^{\fhi}([a,b])$, for every $n \in \N$ sufficiently large. In fact, for every $\la >0$, by using (\ref{max}) and the inequality established in Proposition \ref{proposizione}, we have:
\begin{equation*}
\begin{split}
\I&\bigl[\la\op(f)\bigr]= \integraleab  \fhi \bigl(\la \assolutol \op(f) (x) \assolutor\bigr) \ dx  \le \integraleab \fhi \left(\rapporto  \Vj \assolutol \chi(nx-k) \assolutor \left[n \integralek f(t) \ dt\right]\right) \ dx \\
& \le \integraleab \Vj \varphi \left(2\rapporto \assolutol \chi(nx-k) \assolutor \left[n \integralek f(t) \ dt\right] \right) \ dx \le \varphi\left(2\rapporto \|f\|_{\infty} \|\chi\|_{\infty}\right) \bigl(b-a\bigr) < +\infty,
\end{split}
\end{equation*}
for every $n \in \N$ sufficiently large. Now, let $\ep >0$ be fixed. Thus, for every fixed $\lambda >0$, since the $\fhi$-function $\fhi$ is convex and Theorem \ref{convergenzauniforme} holds, it follows that:
\begin{equation*}
\begin{split}
\I\bigl[\la\bigl(\op(f)-f\bigr)\bigr]&=\integraleab \varphi \bigl(\la \assolutol \op(f)(x)-f(x) \assolutor \bigr) \ dx \le \integraleab \fhi \bigl(\la\|\op(f)-f\|_{\infty}\bigr) \ dx \\
& \le \integraleab \varphi\bigl(\la \ep\bigr) \ dx \le \ep \ \fhi(\la) (b-a),
\end{split}
\end{equation*}
for $n \in \N$ sufficiently large. Thus, by the arbitrariness of $\ep >0$ we get the thesis. 
\newline Now, we can consider the case $\Om =\R$. Let $\gamma$ and $\bar{\gamma}$ be two positive constants such that $supp \ f \subseteq [-\bar{\gamma},\bar{\gamma}]$ and $\gamma > \bar{\gamma}+1$. Thus, whenever $k/n \notin [-\gamma,\gamma]$, $n \in \N$, we have $[k/n,(k+1)/n] \cap [-\bar{\gamma},\bar{\gamma}] = \emptyset$ and therefore:
\[\integralek f(t) \ dt=0.\]
Now, we can prove that $\op(f) \in E^{\fhi}(\R) \subseteq L^{\fhi}(\R)$, for every $n \in \N$. Indeed, for every $\la >0$, we can write what follows: 
\begin{equation*}
\I\bigl[\la \op(f)\bigr]= \integraler \fhi\bigl(\la \assolutol \op(f)(x) \assolutor \bigr) \ dx = \integraler \fhi \left(\la \left\assolutol \dfrac{\displaystyle \Vg \chi(nx-k) n \displaystyle \integralek f(t) \ dt}{\displaystyle \bigvee_{k \in \Z}\chi(nx-k)} \right\assolutor \right) \ dx,
\end{equation*}
where $\J_{n}^{\gamma} := \{k \in \Z: \assolutol k/n \assolutor \le \gamma\}$. Hence, by inequalities (\ref{sup}), and Proposition \ref*{proposizione}, and by using the convexity of $\fhi$ with the fact that $\left|\chi(nx-k)\right| \le \|\chi\|_{\infty}$, for every $k \in \J_{n}^{\gamma}$, we obtain for every $\la >0$:
\begin{equation*}
\begin{split}
\integraler \fhi \bigl(\la \assolutol \op(f)(x) \assolutor \bigr) \ dx &\le \integraler \fhi \left(\rapporto \left \assolutol\Vg \chi(nx-k) n \integralek f(t) \ dt\right \assolutor\right) \ dx \\
& \le \integraler \Vg \fhi \left(2\rapporto  \assolutol \chi(nx-k) \assolutor n \integralek f(t) \ dt\right) \ dx \\
& \le \integraler \Vg \dfrac{\assolutol \chi(nx-k) \assolutor}{\|\chi\|_{\infty}} \fhi \left(2\rapporto\|\chi\|_{\infty} n \integralek f(t) \ dt\right) \ dx \\
& \le \integraler \Vg \dfrac{\assolutol\chi(nx-k)\assolutor}{\|\chi\|_{\infty}} \fhi\left(2\rapporto \|f\|_{\infty} \|\chi\|_{\infty}\right) \ dx \\
& \le \|\chi\|_{\infty}^{-1} \fhi\left(2\rapporto \|f\|_{\infty} \|\chi\|_{\infty}\right) \sum_{k \in \J_{n}^{\gamma}} \integraler \assolutol\chi(nx-k) \assolutor \ dx \\
\end{split}
\end{equation*}
\begin{equation*}
\begin{split}
& = n^{-1}\|\chi\|_{\infty}^{-1} \fhi\left(2\rapporto \|f\|_{\infty} \|\chi\|_{\infty}\right) \sum_{k \in \J_{n}^{\gamma}} \|\chi\|_{1} \\
& =  n^{-1}\|\chi\|_{\infty}^{-1} \fhi\left(2\rapporto \|f\|_{\infty} \|\chi\|_{\infty}\right) \|\chi\|_{1} \cdot L < +\infty, 
\end{split}
\end{equation*} 
where in the last integral, we made the change of variable $y=nx-k$ and $L \le 2n\gamma+1$ represents the number of terms in the above finite sum. 
\newline Now, the goal of the last part of the proof is to prove that: 
\[\lim_{n \rightarrow +\infty} \I\bigl[\la\bigl(\op(f)-f\bigr)\bigr]=\lim_{n \rightarrow +\infty} \integraler \fhi \bigl(\la \assolutol \op(f)(x)-f(x) \assolutor \bigr) \ dx=0,\]
for every $\la >0$, which is equivalent to proving that the family $\bigl(\fhi(\la\assolutol\op(f)-f\assolutor)\bigr)_{n \in \N}$ converges to zero with respect to the $L^1$-norm, for every $\la >0$. In order to establish this, we will exploit the well-known Vitali convergence theorem. Let $\la > 0$ be fixed. First, by Theorem \ref{convergenzauniforme} and the properties of $\fhi$, we have that for every $\ep >0$ there exists $\bar{n} \in \N$ sufficiently large such that: 
\[\fhi\bigl(\la \assolutol \op(f)(x)-f(x) \assolutor\bigr) \le \fhi\bigl(\la \|\op(f)-f\|_{\infty}\bigr) \le \ep,\]
for every $n \ge \bar{n}$ and for all $x \in \R$. Therefore,
\[\assolutol\{x \in \R: \fhi(\la \left|\op(f) (x)-f(x)\right|) > \ep\} \assolutor = \assolutol \emptyset \assolutor =0,\]
for every $n \ge \bar{n}$, where $\assolutol \cdot \assolutor$ denotes the (Lebesgue) measure of the corresponding set. Hence, we proved that the family $\bigl(\fhi(\la \assolutol \op(f)-f \assolutor)\bigr)_{n \in \N}$ converges in measure to zero.
\newline Now, let $\ep >0$ be fixed. Since $\chi$ is absolutely integrable, it is easy to prove that corresponding to 
\[T(\ep):=\dfrac{\ep \|\chi\|_{\infty}}{\fhi\left(2\rapporto\|f\|_{\infty}\|\chi\|_{\infty}\right)(2\gamma+1)},\] there exists a constant $M >0$ (we can assume $M > \bar{\gamma}$), such that
\[n \int_{\assolutol x \assolutor >M}\assolutol \chi(nx-k) \assolutor \ dx< T(\ep),\]
for every $k \in [-\gamma n, \gamma n]$, $n \in \N$. Hence, proceeding as in the previous estimates, we have:
\begin{equation*}
\begin{split}
\integralem \fhi\bigl(\lambda\assolutol \op(f)(x) \assolutor \bigr) \ dx &\le \integralem \Vg \fhi\left(2\rapporto \assolutol\chi(nx-k)\assolutor n \integralek f(t) \ dt \right) \ dx \\
& \le \integralem \Vg \dfrac{\assolutol \chi(nx-k)\assolutor}{\norma\chi\norma_{\infty}} \fhi\left(2\rapporto\norma\chi\norma_{\infty} n \integralek f(t) \ dt\right) \ dx \\
& \le \integralem \Vg \dfrac{\assolutol \chi(nx-k)\assolutor}{\norma\chi\norma_{\infty}} \fhi\left(2\rapporto \norma\chi\norma_{\infty}\norma f\norma_{\infty}\right) \ dx \\
& \le n^{-1}\norma\chi\norma_{\infty}^{-1}\fhi\left(2\rapporto\norma f\norma_{\infty}\norma \chi\norma_{\infty}\right) \sum_{k \in \J_{n}^{\gamma}} n \integralem \assolutol \chi(nx-k) \assolutor \ dx \\
& \le n^{-1}\norma\chi\norma_{\infty}^{-1}\fhi\left(2\rapporto\norma f\norma_{\infty}\norma \chi\norma_{\infty}\right) \sum_{k \in \J_{n}^{\gamma}} T(\ep) \\
& \le \norma\chi\norma_{\infty}^{-1} \fhi\left(2\rapporto\norma f\norma_{\infty}\norma \chi\norma_{\infty}\right) (2\gamma+1) T(\ep)=\varepsilon,
\end{split}
\end{equation*}  
where the last estimate follows from the fact that the sum in the previous inequality consists of at most $2n\gamma+1$ terms, corresponding to the number of intervals $[k/n,(k+1)/n]$, $n \in \N$, having non-empty intersection with $[-\bar{\gamma},\bar{\gamma}]$. Therefore, we proved that for every $\ep >0$ there exists  a set $E_{\ep}:=[-M,M] \subset \R$ for which, for every measurable set $F$, with $ F \cap E_{\ep} = \emptyset$, we get:
\[\int_{F}\fhi\bigl(\la \assolutol\op(f)(x)-f(x)\assolutor\bigr) \ dx= \int_{F}\fhi\bigl(\la \assolutol\op(f)(x)\assolutor\bigr) \ dx \le \integralem \fhi\bigl(\la \assolutol\op(f)(x)\assolutor\bigr) \ dx < \ep,\]
for every $n \in \N$. 
\newline Finally, let $\ep >0$ and $B \subset \R$ be a measurable set such that:
\[\assolutol B \assolutor < \delta(\ep):= \dfrac{\ep}{\fhi\left(2\la \dfrac{\mom}{a_{\chi}}\norma f \norma_{\infty}\right)}.\]
Then, by the convexity of $\fhi$ and observing that \begin{equation}
\label{disugualianzamom}
\mom \ge a_{\chi},
\end{equation}
we obtain:
\begin{equation*}
\integraleb \fhi\bigl(\la \assolutol \op(f)(x)-f(x) \assolutor\bigr) \ dx \le \integraleb \fhi\left(\frac{1}{2}2\la\assolutol\op(f)(x)\assolutor+\frac{1}{2}2\la f(x)\right) \ dx
\end{equation*}
\begin{equation*}
\begin{split}
& \le \frac{1}{2} \integraleb \fhi\bigl(2\la \assolutol\op(f)(x)\assolutor\bigr) \ dx + \frac{1}{2}\integraleb \fhi\bigl(2\la  f(x)\bigr) \ dx \\
& \le \frac{1}{2} \integraleb \fhi\left(2\la \frac{\mom}{a_{\chi}}\norma f \norma_{\infty}\right) \ dx + \frac{1}{2}\integraleb \fhi\bigl(2\la \norma f \norma_{\infty}\bigr) \ dx \\ 
& = \frac{\assolutol B \assolutor}{2} \left[\fhi\left(2\la \frac{\mom}{a_{\chi}}\norma f \norma_{\infty}\right)+ \fhi\bigl(2 \la \norma f \norma_{\infty}\bigr)\right] \\
& \le \assolutol B \assolutor \fhi\left(2\la \dfrac{\mom}{a_{\chi}}\norma f \norma_{\infty}\right) < \ep,
\end{split}
\end{equation*}
for every $n \in \N$.  
Thus, it turns out that the integrals
\[\int_{(\cdot)} \fhi\bigl(\la \assolutol\op(f)(x)-f(x)\assolutor\bigr) \ dx\]
are equi-absolutely continuous. Then, the proof follows by the Vitali convergence theorem, together with the arbitrariness of $\la >0$.    
\end{proof}   
Now, we can finally prove the main theorem of this paper, i.e., a modular convergence theorem for the max-product Kantorovich sampling operators. 
\begin{teorema}
\label{maintheorem}
Let $f \in L^{\fhi}_{+}(\Om)$ be fixed. Then there exists $\la >0$ such that:
\[\lim_{n \rightarrow +\infty} \I\bigl[\la\bigl(\op(f)-f\bigr)\bigr]=0.\] 
\end{teorema} 
\begin{proof}
Let $\ep >0$ be fixed. First of all, we recall that the spaces $C_{+}([a,b])$ and $C_{c}^{+}(\R)$ are modularly dense in $L^{\fhi}_{+}([a,b])$ and $L^{\fhi}_{+}(\R)$, respectively. Thus, corresponding to 
\[T(\ep):=\dfrac{\ep}{2\left(\frac{\norma \chi \norma_{1}}{\mom}+1\right)},\]
there exists $\bar{\la} >0$ and a function $g \in C_{+}([a,b])$, if $\Om=[a,b]$, or $g \in C_{c}^{+}(\R)$, in case of $\Om=\R$, such that:
\begin{equation}
\label{disuguaglianza6}
\I\bigl[\bar{\la}\bigl(f-g\bigr)\bigr] < T(\ep).
\end{equation}
Now, let us choose $\la >0$ such that:
\[\dfrac{\mom}{a_{\chi}}6\la \le \bar{\la}.\]
By the properties of modular $\I$, the convexity of $\fhi$, recalling (\ref{disugualianzamom}) and exploiting Theorem \ref{disuguaglianzamod}, we can write what follows:
\begin{equation*}
\begin{split}
\I\bigl[\la\bigl(\op(f)-f\bigr)\bigr] &\le \frac{1}{3}\bigl\{ \I\bigl[3\la\bigl(\op(f)-\op(g)\bigr)\bigr]+\I\bigl[3\la\bigl(\op(g)-g\bigr)\bigr]+\I\bigl[3\la\bigl(g-f\bigr)\bigr]\bigr\}\\ & \le \dfrac{\norma \chi \norma_{1}}{\mom}\I\left[\dfrac{\mom}{a_{\chi}}6\la\bigl(f-g\bigr)\right] + \I\bigl[3\la\bigl(\op(g)-g\bigr)\bigr]+\I\bigl[3\la\bigl(f-g\bigr)\bigr] \\
& \le \left(\dfrac{\norma \chi \norma_{1}}{\mom}+1\right) \I\bigl[\bar{\la}\bigl(f-g\bigr)\bigr] + \I\bigl[\bar{\la}\bigl(\op(g)-g\bigr)\bigr],
\end{split}
\end{equation*}
for every $n \in \N$ sufficiently large.  
Now, by using (\ref{disuguaglianza6}), we obtain:
\[\I\bigl[\la\bigl(\op(f)-f\bigr)\bigr]\le \dfrac{\ep}{2}+\I\bigl[\bar{\la}\bigl(\op(g)-g\bigr)\bigr],\]
and in view of Theorem \ref{convergenzaforte}, we have for every $n \in \N$ sufficiently large:
\[\I\bigl[\bar{\la}\bigl(\op(g)-g\bigr)\bigr]<\frac{\ep}{2}.\] Thus, the proof follows by the arbitrariness of $\ep >0$.   
\end{proof}
\begin{remark}
As often happens when one deals with max-product type operators, all the convergence results established in this section can be extended to the case of not-necessarily non-negative functions that are bounded from below. Indeed, for any fixed $f: \Om \rightarrow \R$ belonging to suitable spaces, denoting $\inf_{x \in \Om} f(x)=: c$, it turns out that the sequences $(\op(f(\cdot)-c)+c)_{n \in \N}$ converges to $f$ in all the above considered senses.
\end{remark}
\section{Applications to particular cases, examples with special kernels and graphical representations}
\label{esempigrafici}
In this section, we want to show how the convergence results obtained in the previous section can be applied to some special instances of Orlicz spaces. Then, in the last part of the section, we will provide several specific examples of kernels for the max-product Kantorovich sampling operators $(\op)_{n \in \N}$. \newline  Let us begin considering the case that occurs when choosing the convex $\fhi$-function $\fhi(u) := u^{p}$, $u \ge 0$, $1 \le p < +\infty$, which generates the corresponding Orlicz space $L^{\fhi}(\Om)=L^{p}(\Om)$. Indeed, it results $\I[f] = \norma f \norma_{p}^{p}$, i.e., the modular functional $\I$ coincides with the usual $L^{p}$-norm. Further, it is immediate to prove that the $\Delta_{2}$-condition is fulfilled. Hence, it follows that $L^{p}(\Om)$ coincides with the space of its finite elements and the modular convergence and the usual (Luxemburg) norm-convergence are equivalent. From the theory established in the previous section, we can state the following corollaries. 
\begin{cor}
Let $\chi \in L^{1}(\R)$ be a fixed kernel. Then for every non-negative $f$, $g \in L^{p}(\Om)$, $1 \le p < +\infty$, we have:
\[\norma \op(f) -\op(g)\norma_{p} \le \dfrac{2\bigl((\mom)^{p-1}\norma \chi\norma_{1}\bigr)^{1/p}}{a_{\chi}}\cdot \norma f-g \norma_{p},\]
for $n \in \N$ sufficiently large. 
\end{cor}   
\begin{proof}
Since Theorem \ref{disuguaglianzamod} holds with $\I[f]= \norma f \norma_{p}^{p}$, we obtain:
\[\norma\op(f)-\op(g)\norma_{p}^{p} \le 2^{p}a_{\chi}^{-p} \bigl(\mom\bigr)^{p-1} \norma \chi \norma_{1} \norma f -g \norma_{p}^{p},\] 
from which the assertion follows. 
\end{proof}  
In addition, as a direct consequence of Theorem \ref{maintheorem}, we immediately obtain the following $L^{p}$-convergence result.
\begin{cor}
Let $\chi \in L^{1}(\R)$ be a given kernel. For any non-negative function $f \in L^{p}(\Om)$, $1 \le p < +\infty$, there holds:
\[\lim_{n \rightarrow +\infty} \norma \op(f) -f \norma_{p} = 0.\] 
\end{cor}  
Other useful examples of Orlicz spaces are the so-called interpolation spaces $L^{\alpha} \log^{\beta} L(\Om)$. The Zygmund spaces (another name for the interpolation spaces) can be obtained by $\fhi$-functions $\fhi_{\alpha,\beta}(u) := u^{\alpha} \log^{\beta}(u + e)$, $u \ge 0$, for $\alpha \ge 1$ and $\beta > 0$. The corresponding modular functionals
\[I^{\fhi_{\alpha,\beta}}[f]:=\integrale \bigl\assolutol f(x) \bigr\assolutor^{\alpha} \log^{\beta}\bigl(\bigl\assolutol f(x) \bigr\assolutor + e\bigr) \ dx, \quad f \in M(\Om),\] 
are convex, and  $L^{\alpha} \log^{\beta} L(\Om)$ is defined as the set of all (Lebesgue-)measurable functions for which $I^{\fhi_{\alpha,\beta}}[\lambda f] < +\infty$ for some $\lambda >0$, for any $\alpha$ and $\beta$. Applying the modular inequality established in Theorem \ref{disuguaglianzamod} with $\fhi_{\alpha,\beta}(u) = u^{\alpha} \log^{\beta}(u + e)$, we obtain the following corollary, which for simplicity we provide for $\alpha = \beta =1$. 
\begin{cor}
Let $\chi \in L^{1}(\R)$ be a fixed kernel. Then for every non-negative $f$, $g \in L \log L(\Om)$, we have:
\begin{equation*}
\begin{split}
\integrale \bigl\assolutol \op(f)(x) - \op(g)(x) \bigr\assolutor \log\bigl(\lambda \bigl\assolutol \op(f)(x) - \op(g)(x) \bigr\assolutor + e\bigr) \ dx \\ \le \dfrac{2\norma \chi \norma_{1}}{a_{\chi}}\integrale \bigl\assolutol f(x)-g(x) \bigr\assolutor \log\left(\dfrac{\mom}{a_{\chi}}2\la \bigl\assolutol f(x)-g(x) \bigr\assolutor + e\right) \ dx,
\end{split}
\end{equation*}
for $\la >0$ and $n \in \N$ sufficiently large. 
\end{cor}
In addition, we have $L^{\alpha} \log^{\beta} L(\Om)=E^{\fhi_{\alpha,\beta}}(\Om)$, since also $\fhi_{\alpha,\beta}$ satisfies the $\Delta_{2}$-condition. This implies that the modular convergence and norm convergence are equivalent, and the following convergence theorem holds for $\alpha=\beta=1$. 
\begin{cor}
Let $\chi \in L^{1}(\R)$ be a given kernel.  For any non-negative function $f \in L \log L(\Om)$, there holds:
\[\lim_{n \rightarrow +\infty} \integrale \bigl\assolutol \op(f)(x) - f(x) \bigr\assolutor \log\bigl(\la \bigl\assolutol \op(f)(x) - f(x) \bigr\assolutor + e\bigr) \ dx=0,\]
for every $\la >0$ or, equivalently,
\[\lim_{n \rightarrow +\infty}\norma \op(f) -f\norma_{\fhi_{1,1}}=0,\]
where $\norma \cdot \norma_{\fhi_{1,1}}$ is the Luxemburg norm associated to $I^{\fhi_{1,1}}$. 
\end{cor}
As a final example, as mentioned in Section \ref{orlicz}, an important case in which modular and norm convergence are not equivalent occurs at so-called exponential spaces, generated by $\fhi_{\gamma}(u)=e^{u^{\gamma}}-1$, $u \ge 0$, for $\gamma >0$. Indeed, it is easy to prove that in this case the $\Delta_{2}$-condition is not satisfied. Hence, the corresponding space of the finite elements $E^{\fhi_{\gamma}}(\Om)$ is properly contained in $L^{\fhi_{\gamma}}(\Om)$. Now, the modular functional corresponding to $\fhi_{\gamma}$ is given by:
\[I^{\fhi_{\gamma}}[f]:=\integrale \left(e^{\assolutol f(x) \assolutor^{\gamma}}-1\right) \ dx,\]
for any (Lebesgue-)measurable function $f: \Om \rightarrow \R$. As a consequence of Theorem \ref{disuguaglianzamod}, we can establish the following.
\begin{cor}
Let $\chi \in L^{1}(\R)$ be a fixed kernel. Then for every non negative $f$, $g \in L^{\fhi_{\gamma}}(\Om)$, we have:
\begin{equation*}
\begin{split}
\integrale & \biggl(\exp\Bigl(\la \bigl\lvert \op(f) (x) - \op(g)(x) \bigr \rvert\Bigr)^{\gamma}-1\biggr) \ dx \\
& \le \dfrac{\norma \chi \norma_{1}}{\mom} \integrale \biggl(\exp\left(\dfrac{\mom}{a_{\chi}}2\la \bigl \assolutol f(x) - g(x) \bigr \assolutor\right)^{\gamma}-1\biggr) \ dx,
\end{split}
\end{equation*}
for $\la >0$ and $n \in \N$ sufficiently large. 
\end{cor}
In contrast to previous cases, in the latter setting there is no analogous convergence result with respect to the Luxemburg norm, but only a modular convergence theorem.  
\begin{cor}
Let $\chi \in L^{1}(\R)$ be a given kernel. For any non-negative function $f \in L^{\fhi_{\gamma}}(\Om)$, there exists $\la >0$ such that: 
\[\lim_{n \rightarrow +\infty} \integrale \biggl(\exp\Bigl(\la \bigl \assolutol \op(f)(x)-f(x) \bigr \assolutor\Bigr)^{\gamma}-1\biggr) \ dx=0.\]
\end{cor}  
Now, we give some examples of specific kernels for which the results proved in this paper hold. First, we consider the Fejér kernel $\mathcal{F}(x):= \frac{1}{2}$sinc$^{2}(x/2)$, $x \in \R$, (see Figure \ref{fejer}) where the well-known sinc-function is given by:
\[sinc(x):= \begin{cases}
\dfrac{\sin(\pi x)}{\pi x}, & x \in \R \setminus \{0\}, \\
1, & x=0.
\end{cases}\]  
Obviously, $\mathcal{F}$ is a non-negative bounded function, and it belongs to $L^{1}(\R)$, since it turns out that $\mathcal{F}(x) = \mathcal{O}(\assolutol x \assolutor^{-2})$ as $\assolutol x \assolutor \rightarrow +\infty$ (see, e.g., \cite{4.5}). In addition, by Lemma \ref{lemma1}, condition $(\chi_{1})$ holds with $\beta=2$. Then, since $a_{\mathcal{F}}=\mathcal{F}(3/2)>0$, it follows that $(\chi_{2})$ is satisfied, as also $(\chi_{2}')$.
Therefore, for the max-product Kantorovich sampling operators based upon $\mathcal{F}$, all the corollaries of this section hold. 
\begin{figure}[h!]
	\centering
	\includegraphics[width=10cm]{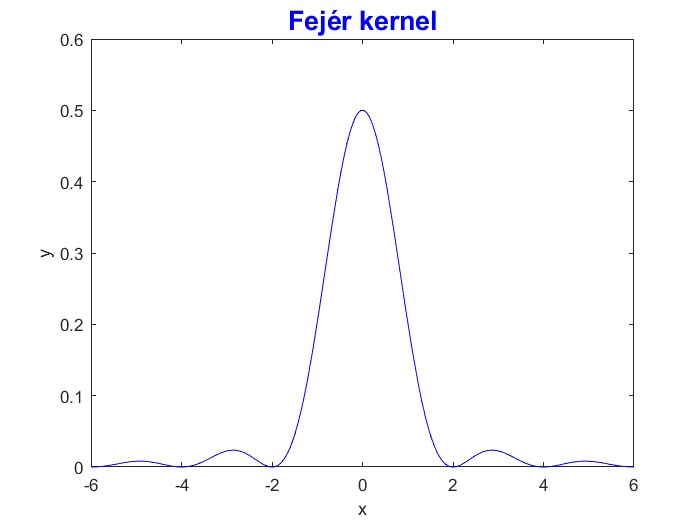}
	\caption{The Fejér kernel $\mathcal{F}$.}
	\label{fejer}
\end{figure} 

The so-called de la Vallée-Poussin kernel $\mathcal{P}$ (see Figure \ref{delavallee}) provides another example of band-limited kernels and it is defined by:
\[\mathcal{P}(x):= \begin{cases}
\dfrac{\sin(x/2)\sin(3x/2)}{9x^{2}/4}, & x\in \R \setminus \{0\}, \\
\dfrac{1}{3}, & x=0.
\end{cases}\]      
We recall that $\mathcal{P}(x) = \mathcal{O}(\assolutol x \assolutor^{-2})$ as $\assolutol x \assolutor \rightarrow +\infty$, which means that $\mathcal{P}$ is absolutely integrable on $\R$ and $m_{\beta}(\mathcal{P}) < +\infty$, i.e., $(\chi_{1})$ holds with $\beta=2$ by Lemma \ref{lemma1}. Then, it is easy to check that $(\chi_{2})$ is satisfied with $a_{\mathcal{P}} \approx 0.1048$, and also $(\chi_{2}')$ is fulfilled (see, e.g., \cite{2}). It follows that the de la Vallée-Poussin kernel $\mathcal{P}$ can be used to define the max-product Kantorovich sampling operators.
\begin{figure}[h!]
	\centering
	\includegraphics[width=10cm]{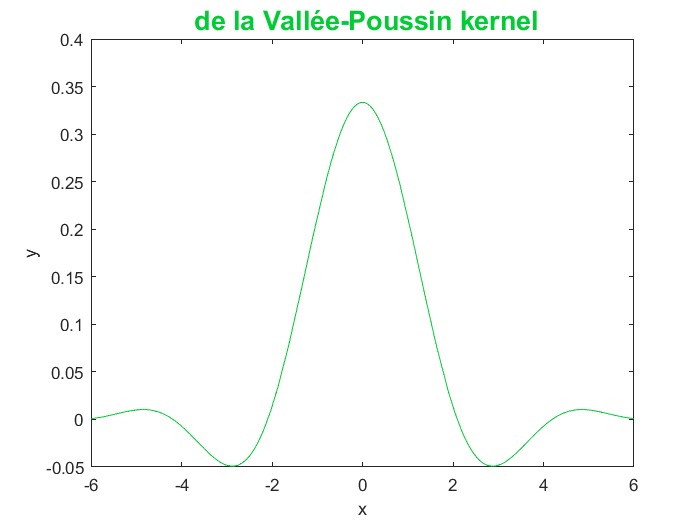}
	\caption{The de la Vallée-Poussin kernel $\mathcal{P}$.}
	\label{delavallee}
\end{figure}
\begin{remark}
\label{remark5.1}
Let $\Om=\R$. We recall that, under suitable assumptions on the kernel $\chi$, the linear version of the Kantorovich sampling operators can be defined as follows:
\begin{equation*}
S_{w}^{\chi}(f)(x)=\sum_{k \in \Z}\chi\left(wx-k\right) \left[w \int_{k/w}^{(k+1)/w} f(t) \ dt \right], \quad x \in \R, \quad w>0,
\end{equation*}
where $f: \R \rightarrow \R$ is a locally integrable function such that the above series is convergent for every $x \in \R$. In various papers (see, e.g., \cite{Bardaro2,4.6}), quantitative estimates for the order of approximation of $(S_{w}^{\chi})_{w>0}$ have been established in case of uniformly continuous and bounded functions, by exploiting the modulus of continuity of the function to be approximated. From the literature, if $f \in C(\R)$, one can define the modulus of continuity of $f$ as follows:
\begin{equation*}
\omega(f, \delta):= \sup \bigl\{\assolutol f(x)-f(y) \assolutor : x,y \in \R, \assolutol x-y \assolutor \le \delta\bigr\}, \quad \delta >0.
\end{equation*}	
In particular, in Theorem 3.2 in \cite{4.6}, the authors proved that, if $\chi=\mathcal{F}$ or $\chi=\mathcal{P}$, then for any $f \in C(\R)$, there exist two suitable positive constants $C_{1}$, $C_{2} >0$ such that:
\begin{equation}
\label{stimalineari}
\norma S_{w}^{\chi} (f) -f \norma_{\infty} \le C_{1} \omega(f, w^{-\beta})+C_{2}w^{-\beta},
\end{equation}
for some $0 < \beta <1$ and $w>0$ sufficiently large. On the other hand, if $f \in C_{+}(\R)$ and $\chi=\mathcal{F}$ or $\chi=\mathcal{P}$, for the max-product version of the same sampling-type operators the following estimate (see Theorem 3.2 in \cite{2}):
\begin{equation*}
\norma \op(f) - f\norma_{\infty} \le \dfrac{2\mom+m_{1}(\chi)}{a_{\chi}}\omega(f, n^{-1}),
\end{equation*}
holds for sufficiently large $n \in \N$. We stress that, by Lemma \ref{lemma1}, $\mom$ and $m_{1}(\chi)$ are both finite. Comparing the above results, it is clear that the Jackson-type estimate for the max-product Kantorovich sampling operators $(\op)_{n \in \N}$ is essentially better than that given in (\ref{stimalineari}) for their linear counterparts.
\end{remark}
Among the kernels with compact support, also known as duration limited kernels, we consider the case generated by the well-known central B-spline of order $n \in \N$, defined by:
\[M_{n}(x):=\dfrac{1}{(n-1)!}\sum_{i=0}^{n}(-1)^{i}\binom{n}{i}\left(\dfrac{n}{2}+x-i\right)_{+}^{n-1}, \quad x \in \R,\]   
where $(x)_{+}:=\max\{x,0\}$ denotes the ``positive part" of $x \in \R$ (see, e.g., \cite{03,8,001}). $M_{n}(x)$ is a non-negative continuous function, and it has support coincident with the interval $[-n/2;n/2]$. This means that, in order to obtain all the convergence results proved in this paper, the central B-spline $M_{n}(x)$ can be used as the kernel of the max-product Kantorovich sampling operators only for $n \ge 4$. Indeed, otherwise $(\chi_{2})$ is not fulfilled. 
\begin{figure}[h!]
\centering
\includegraphics[width=10cm]{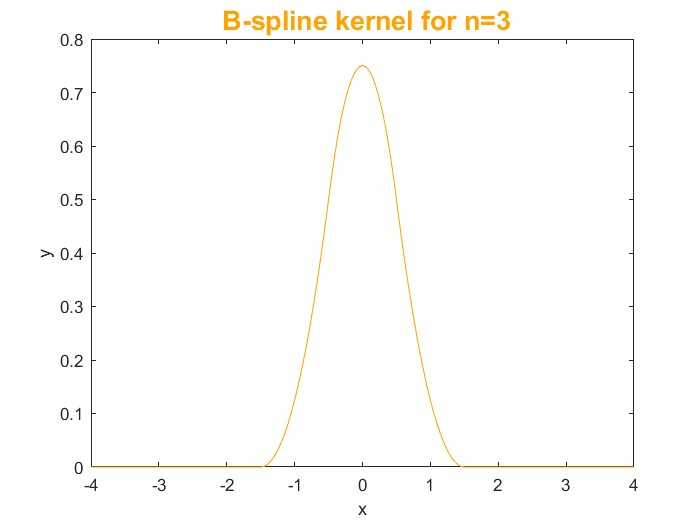}
\caption{The B-spline kernel $M_{3}$.}
\label{bspline}
\end{figure} 
For example, let us now consider the particular case with $n=3$ in detail, for which only $(\chi_{2}')$ is satisfied. For a sake of completeness, it is easy to prove that $M_{3}(x)$ (see Figure \ref{bspline}) can be formally expressed by:
\[M_{3}(x):= \begin{cases}
	\dfrac{3}{4}-x^2, & \assolutol x \assolutor \le \dfrac{1}{2}, \\
	\dfrac{1}{2} \left(\dfrac{3}{2}-\assolutol x \assolutor\right)^{2}, & \dfrac{1}{2} < \assolutol x \assolutor \le \dfrac{3}{2}, \\
	0, & \assolutol x \assolutor > \dfrac{3}{2}. 
\end{cases}\]     
For more details about the above examples and other useful examples of kernels, one can see, e.g., \cite{02,Orlova,Kadak}. 


\section*{Acknowledgments}

{\small The authors are members of the Gruppo Nazionale per l'Analisi Matematica, la Probabilit\`a e le loro Applicazioni (GNAMPA) of the Istituto Nazionale di Alta Matematica (INdAM), of the Gruppo UMI (Unione Matematica Italiana) T.A.A. (Teoria dell'Approssimazione e Applicazioni)}, and of the network RITA (Research ITalian network on Approximation).

\section*{Funding}

{\small The second and the third authors have been partially supported within the (1) 2022 GNAMPA-INdAM Project "Enhancement e segmentazione di immagini mediante operatori di tipo campionamento e metodi variazionali per lo studio di applicazioni biomediche'', (2) "Metodiche di Imaging non invasivo mediante angiografia OCT sequenziale per lo studio delle Retinopatie degenerative dell'Anziano (M.I.R.A.)", funded by the Fondazione Cassa di Risparmio di Perugia (FCRP), 2019 and (3) "National Innovation Ecosystem grant ECS00000041 - VITALITY", funded by the European Union - NextGenerationEU under the Italian Ministry of University and Research (MUR). Moreover, the second author has been partially supported within the 2023 GNAMPA-INdAM Project "Approssimazione costruttiva e astratta mediante operatori di tipo sampling e loro applicazioni".
	
}

\section*{Conflict of interest/Competing interests}

{\small The authors declare that they have no conflict of interest and competing interest.}

\section*{Availability of data and material and Code availability}

{ \small Not applicable.}



\begin{thebibliography}{99}
\bibitem{001}
G. Allasia, R. Cavoretto and A. De Rossi, A class of spline functions for landmark-based image registration, \textit{Math. Meth. Appl. Sci.} \textbf{35} (2012) 923-934.   

\bibitem{004}
C. Bardaro, P. L. Butzer, R. L. Stens and G. Vinti, Kantorovich-type generalized sampling series in the setting of Orlicz spaces, \textit{Sampl. Theory Signal Image Proc.} \textbf{6}(1) (2007) 19-52.    

\bibitem{Bardaro2}
C. Bardaro and I. Mantellini, On convergence properties for a class of Kantorovich discrete operators, \textit{Numer. Funct. Anal. Optim.} \textbf{33}(4) (2012) 374-396.

\bibitem{002} 
C. Bardaro, J. Musielak and G. Vinti, Nonlinear integral operators and applications, \textit{De Gruyter Series in Nonlinear Analysis and Applications}, New York-Berlin \textbf{9} (2003). 

\bibitem{006}
B. Bede, L. Coroianu and S. G. Gal, Approximation and shape preserving properties of the Bernstein operator of max-product kind, \textit{Int. J. Math. Math. Sci.} \textbf{2009} (2009) 590589, doi:10.1155/2009/590589. 

\bibitem{01}
B. Bede, L. Coroianu and S. G. Gal, Approximation by Max-Product Type Operators, (Springer, New York, 2016). 

\bibitem{04}
M. M. Bourke and D. G. Fisher, Solution algorithms for fuzzy relational equations with max-product composition, \textit{Fuzzy sets Syst.} \textbf{94}(1) (1998) 61-69. 

\bibitem{Butzer}
P. L. Butzer, A survey of the Whittaker-Shannon sampling theorem and some of its extensions, \textit{J. Math. Res. Exposition} \textbf{3} (1983) 185-212. 

\bibitem{03}
P. L. Butzer, W. Engels, S. Ries and R. L. Stens, The Shannon sampling series and the reconstruction of signals in terms of linear, quadratic and cubic splines, \textit{SIAM J. Appl. Math.} \textbf{46} (1986) 299-323. 

\bibitem{Butzer4}
P. L. Butzer, A. Fischer and R. L. Stens, Generalized sampling approximation of multivariate signals; general theory, \textit{Atti Sem. Mat. Fis. Univ. Modena} \textbf{41} (1993) 17-37. 

\bibitem{02}
P. L. Butzer and R. J. Nessel, Fourier analysis and approximation I (Academic Press, New York-London, 1971).

\bibitem{Butzer2}
P. L. Butzer, S. Riesz and R. L. Stens, Approximation of continuous and discontinuous functions by generalized sampling series, \textit{J. Approx. Theory} \textbf{50} (1987) 25-39. 

\bibitem{Butzer3}
P. L. Butzer, W. Splettst\"osser and R. L. Stens, The sampling theorem and linear prediction in signal analysis, \textit{Jahresber. Deutsch. Math. Verein.} \textbf{90} (1988) 1-70. 

\bibitem{Butzer5}
P. L. Butzer and R. L. Stens, Reconstruction of signals in $L^{p}(\R)$-space by generalized sampling series based on linear combinations of B-splines, \textit{Integral Trans. Spec. Funct.} \textbf{19}(1) (2008) 35-58. 

\bibitem{cagini}
C. Cagini, D. Costarelli, R. Gujar, M. Lupidi, G. A. Lutty, M. Seracini and G. Vinti, Improvement of retinal OCT angiograms by Sampling Kantorovich algorithms in the assessment of retinal and choroidal perfusion, \textit{Appl. Math. Comput.} \textbf{427}(4) (2022) 127152.  

\bibitem{1}
L. Coroianu, D. Costarelli, S. G. Gal and G. Vinti, The max-product generalized sampling operators: Convergence and quantitative estimates, \textit{Appl. Math. Comput.} \textbf{355} (2019) 173-183.

\bibitem{2}
L. Coroianu, D. Costarelli, S. G. Gal and G. Vinti, Approximation by max-product sampling Kantorovich operators with generalized kernels, \textit{Anal. Appl.} \textbf{19} (2021) 219-244. 

\bibitem{CDGV}
L. Coroianu, D. Costarelli, S. G. Gal and G. Vinti, Connections between the approximation orders of positive linear operators and their max-product counterparts, \textit{Numer. Funct. Anal. Optim.} \textbf{42}(11) (2021) 1263-1286.

\bibitem{2.5}
L. Coroianu and S. G. Gal, Approximation by nonlinear generalized sampling operators of max-product kind, \textit{Sampl. Theory Signal Image Proc.} \textbf{9}(1-3) (2010) 59-75.

\bibitem{3}
L. Coroianu and S. G. Gal, Approximation by max-product sampling operators based on sinc-type kernels, \textit{Sampl. Theory Signal Image Proc.} \textbf{10}(3) (2011) 211-230.

\bibitem{4}
L. Coroianu and S. G. Gal, Saturation results for the truncated max-product sampling operators based on sinc and Fejér-type kernels, \textit{Sampl. Theory Signal Image  Process.} \textbf{11}(1) (2012) 113-132. 

\bibitem{4.4}
L. Coroianu and S. G. Gal, Saturation and inverse results for the Bernstein max-product operator, \textit{Period. Math. Hung.} \textbf{69} (2014) 126-133.

\bibitem{4.5}
L. Coroianu and S. G. Gal, $L^{p}$-approximation by truncated max-product sampling operators of Kantorovich-type based on Fejér kernel, \textit{J. Integr. Equ. Appl.} \textbf{29}(2) (2017) 349-364. 

\bibitem{coroianu1}
L. Coroianu and S. G. Gal, Approximation by truncated max-product operators of Kantorovich-type based on generalized $(\fhi,\psi)$-kernels, \textit{Math. Methods Appl. Sci.} \textbf{41}(17) (2018) 7971-7984. 

\bibitem{4.6}
D. Costarelli, A. M. Minotti and G. Vinti, Approximation of discontinuous signals by sampling Kantorovich series, \textit{J. Math. Anal. Appl.} \textbf{450}(2) (2017) 1083-1103. 

\bibitem{5}
D. Costarelli and R. Spigler, How sharp is the Jensen inequality? \textit{J. Ineq. Appl.} \textbf{2015} (2015) 1-10. 

\bibitem{5.1}
D. Costarelli and G. Vinti, Approximation by max-product neural networks operators of Kantorovich type, \textit{Results Math.} \textbf{69}(1-2) (2016) 505-519. 

\bibitem{5.2}
D. Costarelli and G. Vinti, Convergence for a family of neural network operators in Orlicz spaces, \textit{Math. Nachr.} \textbf{290}(2-3) (2017) 226-235.   

\bibitem{Duman}
O. Duman, Statistical convergence of max-product approximating operators, \textit{Turkish J. Math.} \textbf{34} (2010) 501-514.

\bibitem{6}
S. Y. G\"ung\"or and N. Ispir, Approximation by Bernstein-Chlodowsky operators of max-product kind, \textit{Math. Commun.} \textbf{23} (2018) 205-225.

\bibitem{7}
A. Holhos, Weighted approximation of functions by Meyer-K$\ddot{o}$nig and Zeller operators of max-product type, \textit{Numer. Funct. Anal. Optim.} \textbf{39}(6) (2018) 689-703.

\bibitem{7.8}
A. Holhos, Weighted approximation of functions by Favard operators of max-product type, \textit{Period. Math. Hungar.} \textbf{77}(2) (2018) 340-346. 

\bibitem{Kadak}
U. Kadak, Max-product type multivariate sampling operators and applications to image processing, \textit{Chaos, Solitons and Fractals} \textbf{157} (2022) 111914.

\bibitem{Karakus}
S. Karakus and K. Demirci, Statistical $\sigma$ approximation to max-product operators, \textit{Computer Math. Appl.} \textbf{61} (2011) 1024-1031.

\bibitem{kivinukk}
A. Kivinukk and G. Tamberg, On window methods in generalized Shannon sampling operators, \textit{Appl. Numer. Harmonic Anal.} \textbf{20} (2014) 63-85. 

\bibitem{7.4}
 W. M. Kozlowski, Modular Function Spaces, \textit{Pure Appl. Math.}, Marcel Dekker, New York and Basel (1988).

\bibitem{7.5}
M. A. Krasnosel'ski\v{i} and YA. B. Ruticki\v{i}, Convex functions and Orlicz spaces, P. Noordhoff Ltd., Groningen, The Netherlands, (1961).

\bibitem{7.6}
L. Maligranda, Orlicz spaces and interpolation, IMECC, Campinas, (1989).

\bibitem{7.9}
M. Menekse Yilmaz and G. Uysal, Convergence of singular integral operators in weighted Lebesgue spaces, \textit{Eur. J. Pure Appl. Math.} \textbf{10} no. 2 (2017) 335-347. 

\bibitem{7.1}
J. Musielak, Orlicz spaces and modular spaces, \textit{Lecture notes in Mathematics}, 1034, Springer, Berlin (1983). 

\bibitem{7.3}
J. Musielak and W. Orlicz, On modular spaces, \textit{Studia Math.} \textbf{28} (1959) 49-65. 

\bibitem{Orlova}
O. Orlova and G. Tamberg, On approximation properties of generalized Kantorovich-type sampling operators, \textit{J. Approx. Theory} \textbf{201} (2016) 73-86.

\bibitem{7.7}
M. Rao and Z. Ren, Theory of Orlicz spaces, Dekker Inc., New York-Basel-Hong Kong (1991).

\bibitem{7.2}
M. Rao and Z. Ren, Applications of Orlicz spaces, \textit{Monographs and Textbooks in Pure and Applied Mathematics}, 250, Marcel Dekker Inc., New York (2002). 

\bibitem{ries}
S. Ries and R. L. Stens, Approximation by generalized sampling series, \textit{Constructive theory of functions'84} Sofia (1984) 746-756. 

\bibitem{shannon}
C. Shannon, Communication in the presence of noise, \textit{Proc. I. R. E.} \textbf{37} (1949) 10-21. 

\bibitem{Shen}
RH. Shen and LY. Wei, Convexity of functions produced by Bernstein operators of max-product kind, \textit{Results Math.} \textbf{74} (2019).

\bibitem{8}
M. A. Unser, Ten good reasons for using spline wavelets, in \textit{Optical Science, Engineering and Instrumentation '97} (International Society for Optics and Photonics, 1997) 422-431. 

\bibitem{9}
X. P. Yang, X. G. Zhou and B. Y. Cao, Single-variable term semi-latticized fuzzy relation geometric programming with max-product operator, \textit{Inf. Sci.} \textbf{325} (2015) 271-287.    
\end{thebibliography}
\end{document}